\documentclass[12pt,reqno,draft]{article} 
\usepackage{amsmath,amssymb,amsthm,amsfonts, indentfirst}
\usepackage{enumerate,color,bm}
\usepackage{mathtools}
\makeatletter
\def\@cite#1#2{[{{\bfseries #1}\if@tempswa , #2\fi}]}
\renewcommand{\section}{%
\@startsection{section}{1}{\z@}
{0.5truecm plus -1ex minus -.2ex}%
{1.0ex plus .2ex}{\bfseries\large}}
\def\@seccntformat#1{\csname the#1\endcsname.\ }
\makeatother

\numberwithin{equation}{section} 
\theoremstyle{theorem}
\newtheorem{thm}{Theorem}[section]

\newtheorem{lem}[thm]{Lemma}

\theoremstyle{definition}

\newtheorem{remark}{Remark}[section]

\newtheorem*{prth2.1}{Proof of Theorem 2.1}
\newtheorem*{prth2.2}{Proof of Theorem 2.2}
\newtheorem*{prcor1.2}{Proof of Corollary 1.2}
\newtheorem*{prth1.3}{Proof of Theorem 1.3}
\newtheorem*{op}{Open problem}
\newtheorem*{no}{Note}

\newcommand{\ep}{\varepsilon}
\newcommand{\pa}{\partial}

\newcommand{\R}{\mathbb{R}}

\newcommand{\cl}[1]{{\overline#1}}

\newcommand{\Tmax}{T_{\rm max}}

\newcommand{\red}[1]{{\color{red}#1}}

\begin{document}
\footnote[0]{
    2010{\it Mathematics Subject Classification}\/. 
    Primary: 35B44; 
    Secondary: 35Q92, 92C17.
    }
\begin{center} 
    \Large{{\bf 
    Remarks on finite-time blow-up in a fully parabolic 
    attraction-repulsion chemotaxis system via \\
    reduction to the Keller--Segel system 
    }}%
\end{center}
\vspace{5pt}
\begin{center}
    Yutaro Chiyo, 
    Tomomi Yokota%
      \footnote{Corresponding author.}%
      \footnote{Partially supported by Grant-in-Aid for
      Scientific Research (C), No.\,21K03278.}
    \footnote[0]{
    E-mail: 
    {\tt ycnewssz@gmail.com}, 
    {\tt yokota@rs.tus.ac.jp}
    }\\
    \vspace{12pt}
    Department of Mathematics, 
    Tokyo University of Science\\
    1-3, Kagurazaka, Shinjuku-ku, 
    Tokyo 162-8601, Japan\\
    \vspace{2pt}
\end{center}
\begin{center}    
    \small \today
\end{center}

\vspace{2pt}
\newenvironment{summary}
{\vspace{.5\baselineskip}\begin{list}{}{%
     \setlength{\baselineskip}{0.85\baselineskip}
     \setlength{\topsep}{0pt}
     \setlength{\leftmargin}{12mm}
     \setlength{\rightmargin}{12mm}
     \setlength{\listparindent}{0mm}
     \setlength{\itemindent}{\listparindent}
     \setlength{\parsep}{0pt}
     \item\relax}}{\end{list}\vspace{.5\baselineskip}}
\begin{summary}
{\footnotesize {\bf Abstract.} 
This paper deals with the fully parabolic attraction-repulsion 
chemotaxis system
%
\begin{align*}
\begin{cases}
      u_t=\Delta u-\chi\nabla \cdot (u\nabla v)+\xi \nabla\cdot(u \nabla w),
            & x \in \Omega,\ t>0,
    \\
      v_t=\Delta v-v+u,
            & x \in \Omega,\ t>0,
    \\
      w_t=\Delta w-w+u,
            & x \in \Omega,\ t>0
\end{cases}
\end{align*}
%
under homogeneous Neumann boundary conditions and initial conditions, 
where $\Omega$ is an open ball in $\R^n$ ($n \ge 3$), 
$\chi, \xi>0$ are constants. 
When $w=0$, finite-time blow-up in the 
corresponding Keller--Segel system has already been obtained. 
However, finite-time blow-up in the above attraction-repulsion 
chemotaxis system has not yet been established 
except for the case $n=3$. 
This paper provides an answer to this open problem by using
a transformation which leads to a system presenting structural advantages respect to the original.}
\end{summary}
\vspace{10pt}

\newpage

\section{Introduction} \label{Sec1}

In this paper we consider the fully parabolic attraction-repulsion 
chemotaxis system 
%
\begin{align}\label{P}
\begin{cases}
      u_t=\Delta u-\chi\nabla \cdot (u\nabla v)+\xi \nabla\cdot(u \nabla w),
            & x \in \Omega,\ t>0,
    \\[1.05mm]
      v_t=\Delta v-v+u,
            & x \in \Omega,\ t>0,
    \\[1.05mm]
      w_t=\Delta w-w+u,
            & x \in \Omega,\ t>0,
    \\[1.5mm]
      \nabla u \cdot \nu=\nabla v \cdot \nu=\nabla w \cdot \nu=0,
            & x \in \pa\Omega,\ t>0,
    \\[1.05mm]
      u(x,0)=u_0(x),\ v(x,0)=v_0(x),\ w(x,0)=w_0(x),
            & x \in \Omega,
\end{cases}
\end{align}
%
where $\Omega \coloneqq B(0, R) \subset \R^n$ ($n \ge 3$) is an open ball 
centered at the origin with radius $R>0$; 
$\chi, \xi>0$ are constants; 
$\nu$ is the outward normal vector to $\pa \Omega$. 
Moreover, the initial data $u_0, v_0, w_0$ are supposed to be 
radially symmetric and positive functions 
which satisfy that 
$$(u_0, v_0, w_0) 
     \in C^0(\cl{\Omega}) \times W^{1, \infty}(\Omega) 
          \times W^{1, \infty}(\Omega).$$
%
%

%
In biology, the functions $u, v$ and $w$ represent the cell density, 
the concentration of attractive and repulsive chemical substances, respectively. 
The system \eqref{P} is one of variations of the chemotaxis system 
proposed by Keller and Segel \cite{KS-1970} (see also Hillen--Painter \cite{HP-2009}). 
The first equation in \eqref{P} describes 
the time evolution of the cell density
in response to its own chemical attractants and repellents. 
Specifically, it implies that the cell movement 
is directed toward a higher concentration of the attractive signal and 
away from the repulsive signal.
\vspace{2mm}

On the other hand, in mathematics, it is important to consider 
whether a solution of the system \eqref{P} can blow up or not. 
In this paper we show finite-time blow-up of a solution to 
the system \eqref{P}. 
The classical parabolic--elliptic Keller--Segel 
system and the parabolic--elliptic--elliptic attraction-repulsion chemotaxis 
system have been investigated in many literatures 
on chemotaxis systems (see e.g., Arumugam--Tyagi \cite{AT}). 
Before presenting the main result, we give an overview of 
known results about some problems related to \eqref{P}.
\vspace{2mm}

We first focus on the parabolic--parabolic chemotaxis system
%
\begin{align}\label{pp}
\begin{cases}
      u_t=\Delta u-\nabla \cdot (u\nabla v),
    \\[1.05mm]
      v_t=\Delta v-v+u.
\end{cases}
\end{align}
%
The first result on unboundedness of solutions to \eqref{pp} 
was established by Winkler \cite{W-2010-1}. 
After that, Winkler \cite{W-2013} succeeded in showing that 
a solution of \eqref{pp} blows up in finite time  
under the condition that for all $m>0$, $A>0$, 
the initial data belongs to 
$\mathcal{B}(m, A)$ which is the set of 
radially symmetric positive functions 
$(\varphi, \psi) \in C^0(\cl{\Omega}) \times W^{1, \infty}(\Omega)$ 
satisfying $\int_\Omega \varphi=m$, 
$\|\psi\|_{W^{1, 2}(\Omega)} \le A$ and 
$\mathcal{F}(\varphi, \psi) \le -K(m, A)$ with some $K(m, A)>0$, 
where $\mathcal{F}$ is the energy functional. 
Moreover, it was shown in \cite{W-2013} 
that $\mathcal{B}(m, A)$ is dense in the space of all 
radially symmetric positive functions in 
$C^0(\cl{\Omega}) \times W^{1, \infty}(\Omega)$ 
with respect to the topology in $L^p(\Omega) \times W^{1,2}(\Omega)$ 
for all $p \in (1, \frac{2n}{n+2})$. 
Also, some related works which derive lower bound of blow-up time 
can be found in \cite{FMV-2015, MVV-2015, MVV-2016}. 
\vspace{2mm}

Secondly, we turn our eyes into the parabolic--elliptic 
chemotaxis system with signal-dependent sensitivity, 
%
\begin{align}\label{pe}
\begin{cases}
      u_t=\Delta u-\nabla \cdot (u\nabla v),
    \\[1.05mm]
      0=\Delta v-v+u.
\end{cases}
\end{align}
%
As to this system, Nagai \cite{N-1995} derived that 
a radially symmetric solution blows up in finite time 
under some condition for the energy function and the moment of $u$ 
in two or more space dimensions. 
After that, in the two-dimensional setting, 
Nagai \cite{N-2001} proved that if $\int_\Omega u_0(x)|x-x_0|^2\,dx$ 
for $x_0 \in \Omega$ is sufficiently small and 
$\int_\Omega u_0(x)\,dx>\frac{8\pi}{\chi}$ holds, 
then there exists a non-radial solution which blows up in finite time. 

\vspace{2mm}

We now shift our attention to the parabolic--elliptic--elliptic 
version of the attraction-repulsion chemotaxis system
%
\begin{align}\label{pee}
\begin{cases}
      u_t=\Delta u-\chi\nabla \cdot (u\nabla v)
            +\xi \nabla\cdot(u \nabla w),
    \\[1.05mm]
      0=\Delta v+\alpha u-\beta v,
    \\[1.05mm]
      0=\Delta w+\gamma u-\delta w,
\end{cases}
\end{align}
%
where $\chi, \xi, \alpha, \beta, \gamma, \delta>0$ are constants. 
Existence of a solution which blows up in finite time 
was studied by \cite{LL-2016, TW-2013, YGZ-2017}. 
More precisely, in the two-dimensional setting, 
Tao and Wang \cite{TW-2013} derived finite-time blow-up 
under the conditions that $\int_\Omega u_0(x)|x-x_0|^2\,dx$ 
for $x_0 \in \Omega$ is sufficiently small and that
%
\begin{itemize}
\item[] \hspace{-8mm} (i)
     $\chi\alpha-\xi\gamma>0$, \quad $\delta=\beta$\quad and 
     \quad $\int_\Omega u_0(x)\,dx>\frac{8\pi}{\chi\alpha-\xi\gamma}$.
\end{itemize}
%
The idea in \cite{TW-2013} is to reduce \eqref{pee} to the form \eqref{pe} 
by introducing the linear combination $z \coloneqq \chi v-\xi w$. 
Also, in the two-dimensional setting, Li and Li \cite{LL-2016} 
extended the above (i) to the following two conditions:
%
\begin{itemize}
\item[] \hspace{-8mm} (ii)
     $\chi\alpha-\xi\gamma>0$, \quad $\delta \ge \beta$\quad and 
     \quad $\int_\Omega u_0(x)\,dx>\frac{8\pi}{\chi\alpha-\xi\gamma}$;
\item[] \hspace{-8mm} (iii)
     $\chi\alpha\delta-\xi\gamma\beta>0$, \quad $\delta<\beta$ \quad 
     and \quad 
     $\int_\Omega u_0(x)\,dx>\frac{8\pi}{\chi\alpha\delta-\xi\gamma\beta}$.
\end{itemize}
%
After that, Yu, Guo and Zheng \cite{YGZ-2017} 
improved the above (iii) by replacing 
$\chi\alpha\delta-\xi\gamma\beta$ with $\chi\alpha-\xi\gamma$ 
in the first and third conditions in (iii) 
and filled the gap between the above (ii) and (iii); 
note that 
any relationship between 
$\delta$ and $\beta$ is no longer necessary. 
On the other hand, in the two dimensional setting, 
Viglialoro \cite{V-2019} provided 
an explicit lower bound of blow-up time for the system \eqref{pee}. 
\vspace{2mm}

For the fully parabolic attraction-repulsion system with positive parameters $\alpha, \beta, \gamma, \delta$, i.e., 
the fully parabolic version of \eqref{pee}, 
Lankeit~\cite{L-pre} succeeded in establishing existence of radially symmetric solutions blowing up 
at some finite time under the condition that $\chi\alpha-\xi\gamma>0$ 
without any restriction on $\beta, \delta$ in the three-dimensional setting.

\vspace{2mm}

%
In summary, finite-time blow-up has been shown 
for the {\it parabolic--elliptic--elliptic} attraction-repulsion chemotaxis 
system \eqref{pee}. 
However, finite-time blow-up in the {\it fully parabolic} system \eqref{P} 
has not been obtained yet except for the case $n=3$. 
%
\vspace{2mm}

The purpose of this paper is to give an answer to the above open problem, 
that is, to establish finite-time blow-up 
in the {\it fully parabolic} system \eqref{P}. 
The strategy for proving finite-time blow-up is 
to apply the method in \cite{W-2013} to 
the system \eqref{P} via the linear combination of the solution components $v,w$ 
such that 
$z \coloneqq \chi v-\xi w$. 
After this transformation, the analysis reduces to citing a well-known result 
from the literature, but a new information about blow-up in \eqref{P} is 
obtained. 


\section{Main results and their proofs}
In this section we give two main theorems. 
The first one asserts finite-time blow-up in the system \eqref{P}. 
The statement and proof read as follows.
%
\begin{thm}\label{thm1}
Let\/ $\Omega \coloneqq B(0, R) \subset \R^n$ $(n \ge 3)$ 
be an open ball centered at the origin with radius $R>0$. 
Let $m>0$, $A>0$ and $\chi>\xi$. 
Then there exist constants $T=T(m, A)>0$ and $K=K(m, A)>0$ such that 
if $(u_0, v_0, w_0)$ belongs to the set 
%
\begin{align}\label{BmA}
    &\mathcal{C}(m,A)
      :=\Big\{
         (u_0, v_0, w_0) 
             \in C^0(\cl{\Omega}) \times W^{1, \infty}(\Omega) 
                  \times W^{1, \infty}(\Omega)
         \,\Big|\, \notag\\
    &\hspace{2.5cm}
         u_0\ {\it and}\ \chi v_0-\xi w_0\ 
         {\it are\ radially\ symmetric\ and\ positive\ in}\ \cl{\Omega}
         \notag\\
    &\hspace{2.5cm}
         {\it with}\ \int_\Omega u_0=m,\ 
         \|\chi v_0-\xi w_0\|_{W^{1, 2}(\Omega)} \le A\ 
         {\it and}\ 
         \mathcal{G}(u_0, v_0, w_0) \le -K
      \Big\},
\end{align}
%
where $\mathcal{G}$ is the energy functional defined as 
%
\begin{align} \label{ene}
      \mathcal{G}(u_0, v_0, w_0)
      &:=\frac{1}{2}\int_\Omega |\nabla (\chi v_0-\xi w_0)|^2
         +\frac{1}{2}\int_\Omega (\chi v_0-\xi w_0)^2
      \notag\\
      &\qquad\quad
         -(\chi-\xi)\int_\Omega u_0(\chi v_0-\xi w_0)
         +(\chi-\xi)\int_\Omega u_0\ln u_0, 
\end{align}
%
then the corresponding solution $(u, v, w)$ of the system \eqref{P} 
blows up before or at time $T$.
\end{thm}
%

\begin{remark}
For $n = 3$, the above result is covered 
by~\cite[Theorem~1.1 with $\beta=\delta=1$]{L-pre}. 
\end{remark}

Before proving the above theorem,  
we give the local 
solvability to clarify blow-up. 
\begin{lem}[{\cite[Lemma~3.1]{TW-2013}}] \label{LSE}
Let\/ $(u_0, v_0, w_0) 
             \in C^0(\cl{\Omega}) \times W^{1, \infty}(\Omega) 
                  \times W^{1, \infty}(\Omega)$. 
 Then there exists $\Tmax \in (0,\infty]$ such that 
 \eqref{P} possesses a unique classical solution 
 $(u, v, w)$ such that
%
\begin{align*}
    &u,v,w\in C^0(\overline{\Omega}\times[0,\Tmax))
           \cap C^{2,1}(\overline{\Omega}\times(0,\Tmax)),
\end{align*}
%
 and
%
    \begin{align*}
            {\it if}\ \Tmax<\infty,
    \quad 
            {\it then}\ \lim_{t \nearrow \Tmax} 
                           \|u(\cdot,t)\|_{L^\infty(\Omega)}=\infty.
    \end{align*}
%
In particular, if $u_0 \ge 0$, then 
$u(\cdot, t)\ge0$ for all $t \in (0, \Tmax)$. 
Moreover, if $u_0, v_0, w_0$ are radially symmetric in $\cl{\Omega}$, 
then so are $u, v, w$.
\end{lem}

%
\begin{prth2.1} 
We introduce the linear combination of the variables $(v,w)$ such that 
$z \coloneqq \chi v-\xi w$. Then the system \eqref{P} is rewritten as 
\begin{align}\label{ReP}
\begin{cases}
      u_t=\Delta u-\nabla \cdot (u\nabla z),
            & x \in \Omega,\ t \in (0, \Tmax),
    \\[1.05mm]
      z_t=\Delta z-z+(\chi-\xi)u,
            & x \in \Omega,\ t \in (0, \Tmax),
    \\[1.5mm]
      \nabla u \cdot \nu=\nabla z \cdot \nu=0,
            & x \in \pa\Omega,\ t \in (0, \Tmax),
    \\[1.05mm]
      u(x,0)=u_0(x),\ z(x,0)=z_0(x),
            & x \in \Omega,
\end{cases}
\end{align}
where $z_0 \coloneqq \chi v_0-\xi w_0$ and $\chi>\xi$.  
We also define the energy functional for \eqref{ReP} as
\begin{align*}
    \mathcal{F}(u,z)
    :=
      \frac{1}{2}\int_\Omega |\nabla z|^2+\frac{1}{2}\int_\Omega z^2
      -(\chi-\xi)\int_\Omega uz+(\chi-\xi)\int_\Omega u\ln u.
\end{align*}
Then the functional $\mathcal{F}$ satisfies the energy inequality
\begin{align}\label{ei}
      \frac{d}{dt}\mathcal{F}(u(\cdot, t), z(\cdot, t))
    \le
      -\mathcal{D}(u(\cdot, t), z(\cdot, t))
    \quad
      {\rm for\ all}\ t \in (0, \Tmax), 
\end{align}
where $\mathcal{D}$ is the dissipation rate defined as
\begin{align*}
      \mathcal{D}(u, z)
    \coloneqq \int_\Omega z_t^2
      +(\chi-\xi)\int_\Omega u \cdot \big|\nabla \ln u-\nabla z\big|^2
\end{align*}
and $\Tmax \in (0, \infty]$ is the maximal existence time of 
the solution $(u, z)$ to the system \eqref{ReP}.
Indeed, in order to confirm \eqref{ei}, multiplying the second equation in \eqref{ReP} 
by $z_t$, we have
\[
\int_\Omega {z_t}^2+ \frac{d}{dt}\Big[
       \frac{1}{2}\int_\Omega |\nabla z|^2+\frac{1}{2}\int_\Omega z^2
       -(\chi-\xi)\int_\Omega u(z-\ln u)\Big] = (\chi-\xi)\int_\Omega u_t(\ln u-z).
\]
This together with 
\[
\int_\Omega u_t(\ln u-z)
=-\int_\Omega u\cdot \big|\nabla \ln u-\nabla z\big|^2 
\]
implies \eqref{ei}.
Now we shall verify the conditions for blow-up in \cite[Theorem 1.1]{W-2013} which asserts 
that for all $m>0$ and $A>0$ there exist $T=T(m, A)>0$ and 
$K=K(m, A)>0$ 
such that if $(u_0, z_0)$ belongs to the set 
%
\begin{align}\label{BmA2}
    &\mathcal{B}(m,A)
      :=\Big\{
         (u_0, z_0) 
             \in C^0(\cl{\Omega}) \times W^{1, \infty}(\Omega)
         \,\Big|\, \notag\\
    &\hspace{2.5cm}
         u_0, z_0\ 
         {\rm are\ radially\ symmetric\ and\ positive\ in}\ \cl{\Omega}
         \notag\\
    &\hspace{2.5cm}
         {\rm with}\ \int_\Omega u_0=m,\ 
         \|z_0\|_{W^{1, 2}(\Omega)} \le A\ 
         {\rm and}\ 
         \mathcal{F}(u_0, z_0) \le -K
      \Big\},
\end{align}
%
then the corresponding solution $(u, z)$ of \eqref{ReP} blows up 
before or at time $T$, provided that $\chi-\xi=1$.  
Recalling that $z=\chi v-\xi w$, $z_0=\chi v_0-\xi w_0$ 
and 
\[
\mathcal{F}(u_0, z_0)=\mathcal{G}(u_0, v_0, w_0),
\] 
we see by the assumption of Theorem \ref{thm1} 
that all the above conditions for blow-up are satisfied 
at least in the case $\chi-\xi=1$. In the case  $\chi-\xi>0$, 
we can develop the proof similarly. 
Indeed, as in~\cite[Theorem~5.1]{W-2013},  
we first find $c_1=c_1(m, A, n)>0$ such that
\begin{equation}\label{F-D}
   \mathcal{F}(\tilde{u},\tilde{z}) 
   \ge -c_1 \left(\mathcal{D}^\theta(\tilde{u},\tilde{z})+1\right) 
\end{equation}
with $\theta:=\frac{n+2}{n+4}$ for all $(\tilde{u},\tilde{z}) \in 
\{(u,z) \in C^1(\cl{\Omega})\times C^2(\cl{\Omega}) \mid 
u,z$ {\rm are radially symmetric and positive with} 
$\nabla z\cdot\nu=0$ {\rm on} $\partial\Omega$ {\rm and} 
$\int_\Omega u=m$, $\int_\Omega z \le M$, $z(x) \le B|x|^{-\kappa}$ 
{\rm for all} $M, B>0$ {\rm and} $\kappa>n-2\}$. 
In particular, \eqref{F-D} 
can be applied to $(\tilde{u},\tilde{z}):=(u(\cdot,t), z(\cdot,t))$ for each $t \in (0,\Tmax)$ 
with some suitable $M, B>0$. Setting 
\[
     y(t):=-\mathcal{F}(u(\cdot,t), z(\cdot,t)), \quad t \in [0,\Tmax),
\]
we next combine \eqref{F-D} with \eqref{ei} and consequently we have
\[
   y'(t) \ge c_2 y^\frac{1}{\theta}(t) \quad \text{for all}\ t\in (0,\Tmax)
\]
with some $c_2=c_2(m, A, n)>0$. This yields that
\[
   y(t) \ge y(0)\cdot\left(1-\frac{1-\theta}{\theta}c_2y^{\frac{1-\theta}{\theta}}(0)\,t\right)^{-\frac{\theta}{1-\theta}} \quad \text{for all}\ t\in (0,\Tmax).
\]
In particular, this leads us to the conclusion $\Tmax<\infty$. \qed
\end{prth2.1}
%

We next give and show the second main theorem, 
which guarantees that if the set $\mathcal{C}(m, A)$ 
defined in \eqref{BmA} equipped 
with a suitable topology, then it is dense in the space of 
radially symmetric positive functions. 
%
\begin{thm}\label{thm2}
Let\/ $\Omega \coloneqq B(0, R) \subset \R^n$ $(n \ge 3)$ 
be an open ball centered at the origin with radius $R>0$. 
Let $p \in (1, \frac{2n}{n+2})$. 
Then for all $m>0$, $A>0$, the set $\red{\mathcal{C}}(m, A)$ 
defined in \eqref{BmA} is dense in the space
%
\begin{align*}
    &Y := \Big\{
         (u, v, w) \in C^0(\cl{\Omega}) \times W^{1, \infty}(\Omega) 
                          \times W^{1, \infty}(\Omega) 
         \,\Big|\,\notag\\
    &\hspace{1.35cm}
         u\ {\it and}\ \chi v -\xi w\ 
         {\it are\ radially\ symmetric\ and\ positive\ in}\ \cl{\Omega}
      \Big\}
\end{align*}
%
with respect to the topology in 
$L^p(\Omega) \times W^{1,2}(\Omega) \times W^{1,2}(\Omega)$, 
that is, for all $(u_0, v_0, w_0) \in \mathcal{C}(m, A)$ 
and all $\ep>0$ there exists $(u_{0\ep}, v_{0\ep}, w_{0\ep}) \in Y$ 
such that 
%
\begin{align*}
      \|u_{0\ep}-u_0\|_{L^p(\Omega)}
    +\|v_{0\ep}-v_0\|_{W^{1,2}(\Omega)}
    +\|w_{0\ep}-w_0\|_{W^{1,2}(\Omega)}<\ep.
\end{align*}
%
In particular, the solution $(u_\ep, v_\ep, w_\ep)$ of the system \eqref{P} 
with initial data $(u_\ep, v_\ep, w_\ep)|_{t=0}=(u_{0\ep}, v_{0\ep}, w_{0\ep})$ 
blows up in finite time.
\end{thm}
%

%
\begin{prth2.2}
Let $n \ge 3$, $p \in (1, \frac{2n}{n+2})$, $m>0$ and $A>0$. 
Then we see from \cite[Theorem 1.2]{W-2013} that 
the set $\mathcal{B}(m, A)$ defined in \eqref{BmA2} is dense 
in the space
%
\begin{align*}
    &X:=\Big\{
         (u, z) \in C^0(\cl{\Omega}) \times W^{1, \infty}(\Omega)  
         \,\Big|\,
         u, z\ {\rm are\ radially\ symmetric\ and\ positive\ in}\ 
      \cl{\Omega}\Big\},
\end{align*}
%
and that for all $(u_0, z_0) \in X$ and all $\ep>0$ 
there exists $(u_{0\ep}, z_{0\ep}) \in X$ such that 
the solution $(u_\ep, z_\ep)$ of the system \eqref{ReP} 
with initial data $(u_\ep, z_\ep)|_{t=0}=(u_{0\ep}, z_{0\ep})$ 
blows up in finite time, which concludes the proof. \qed
\end{prth2.2}
%

%

\begin{op} 
Consider the fully parabolic attraction-repulsion chemotaxis system with 
positive parameters $\alpha, \beta, \gamma, \delta$: 
\begin{align*}
\begin{cases}
      u_t=\Delta u-\chi\nabla \cdot (u\nabla v)+\xi \nabla\cdot(u \nabla w),
            & x \in \Omega,\ t>0,
    \\[1.05mm]
      v_t=\Delta v-\beta v+ \alpha u,
            & x \in \Omega,\ t>0,
    \\[1.05mm]
      w_t=\Delta w-\delta w+\gamma u,
            & x \in \Omega,\ t>0.
\end{cases}
\end{align*}
This cannot be reduced to the Keller--Segel system as in \eqref{ReP} 
via the transformation $z=\chi u -\xi w$ in the case 
\[
 \beta \neq \delta.
\]
In this case, Lankeit~\cite{L-pre} established 
finite-time blow-up in the three dimensional-setting. 
However, for $n \ge 4$, finite-time blow-up is left as an open problem.
\end{op}

\begin{no}
After the completion of this paper, 
we confirmed that 
Fujie and Suzuki \cite[Remark 1.5]{FS-2019} derived
Theorem~\ref{thm1} by the same 
transformation $z=\chi u -\xi w$. 
\end{no}


\begin{thebibliography}{10}

\bibitem{AT}
G.~Arumugam and J.~Tyagi.
\newblock {K}eller--{S}egel chemotaxis models: A review.
\newblock {\em Acta Appl.\ Math.}, {\bf 171}(6):82pp., 2021.

\bibitem{FMV-2015}
M.~A. Farina, M.~Marras, and G.~Viglialoro.
\newblock On explicit lower bounds and blow-up times in a model of chemotaxis.
\newblock {\em Discrete Contin.\ Dyn.\ Syst.}, (Dynamical systems, differential
  equations and applications.\ 10th AIMS Conference.\ Suppl.):409--417, 2015.

\bibitem{FS-2019}
K.~Fujie, T.~Suzuki. 
\newblock Global existence and boundedness in a fully parabolic 2D 
attraction-repulsion system: chemotaxis--dominant case, 
\newblock {\em Adv.\ Math.\ Sci.\ Appl.}, {\bf 28}:1--9, 2019.

\bibitem{HP-2009}
T.~Hillen and K.~J. Painter.
\newblock A user's guide to {PDE} models for chemotaxis.
\newblock {\em J. Math.\ Biol.}, {\bf 58}(1--2):183--217, 2009.

\bibitem{KS-1970}
E.~F. Keller and L.~A. Segel.
\newblock Initiation of slime mold aggregation viewed as an instability.
\newblock {\em J. Theoret.\ Biol.}, {\bf 26}(3):399--415, 1970.

\bibitem{L-pre}
J.~Lankeit.
\newblock Finite-time blow-up in the three-dimensional fully parabolic
  attraction-dominated attraction-repulsion chemotaxis system.
\newblock {\em arXiv: 2103.17044 [math.AP]}, 2021.

\bibitem{LL-2016}
Y.~Li and Y.~Li.
\newblock Blow-up of nonradial solutions to attraction-repulsion chemotaxis
  system in two dimensions.
\newblock {\em Nonlinear Anal.\ Real World Appl.}, {\bf 30}:170--183, 2016.

\bibitem{MVV-2015}
M.~Marras, S.~Vernier-Piro, and G.~Viglialoro.
\newblock Lower bounds for blow-up in a parabolic--parabolic {K}eller--{S}egel
  system.
\newblock {\em Discrete Contin.\ Dyn.\ Syst.}, (Dynamical systems, differential
  equations and applications.\ 10th AIMS Conference.\ Suppl.):809--816, 2015.

\bibitem{MVV-2016}
M.~Marras, S.~Vernier-Piro, and G.~Viglialoro.
\newblock Blow-up phenomena in chemotaxis systems with a source term.
\newblock {\em Math.\ Methods Appl.\ Sci.}, {\bf 39}(11):2787--2798, 2016.

\bibitem{N-1995}
T.~Nagai.
\newblock Blow-up of radially symmetric solutions to a chemotaxis system.
\newblock {\em Adv.\ Math.\ Sci.\ Appl.}, {\bf 5}(2):581--601, 1995.

\bibitem{N-2001}
T.~Nagai.
\newblock Blowup of nonradial solutions to parabolic--elliptic systems modeling
  chemotaxis in two-dimensional domains.
\newblock {\em J. Inequal.\ Appl.}, {\bf 6}(1):37--55, 2001.

\bibitem{TW-2013}
Y.~Tao and Z-A. Wang.
\newblock Competing effects of attraction vs.\ repulsion in chemotaxis.
\newblock {\em Math.\ Models Methods Appl.\ Sci.}, {\bf 23}(1):1--36, 2013.

\bibitem{V-2019}
G.~Viglialoro.
\newblock Explicit lower bound of blow-up time for an attraction-repulsion
  chemotaxis system.
\newblock {\em J. Math.\ Anal.\ Appl.}, {\bf 479}(1):1069--1077, 2019.

\bibitem{W-2010-1}
M.~Winkler.
\newblock Aggregation vs.\ global diffusive behavior in the higher-dimensional
  {K}eller--{S}egel model.
\newblock {\em J. Differential Equations}, {\bf 248}(12):2889--2905, 2010.

\bibitem{W-2013}
M.~Winkler.
\newblock Finite-time blow-up in the higher-dimensional parabolic--parabolic
  {K}eller--{S}egel system.
\newblock {\em J. Math.\ Pures Appl.\ $(9)$}, {\bf 100}(5):748--767, 2013.

\bibitem{YGZ-2017}
H.~Yu, Q.~Guo, and S.~Zheng.
\newblock Finite time blow-up of nonradial solutions in an attraction-repulsion
  chemotaxis system.
\newblock {\em Nonlinear Anal.\ Real World Appl.}, {\bf 34}:335--342, 2017.

\end{thebibliography}
\end{document}